\documentclass[12pt]{article}

\usepackage{amsmath}       
\usepackage{amssymb}       
\usepackage{chapterbib}    
\usepackage{epsfig}        
\usepackage{color}
\usepackage{comment}

\usepackage{bm,bbm}

\newcommand{\EOD}{\end{document}}

\usepackage{color}

\def\Frac{\displaystyle \frac}
\def\Int{ \displaystyle \int}
\def\Sum{ \displaystyle \sum}

\def\div{{ \rm div}}

\def\curl{{\rm curl}}

\def\DIV{{\cal DIV}}
\def\DIVt{\widetilde{\cal DIV}}
\def\GRAD{{\cal GRAD}}
\def\GRADt{\widetilde{\cal GRAD}}
\def\CURL{{\cal CURL}}
\def\CURLt{\widetilde{\cal CURL}}

\def\p{\partial}

\def\dx{\,{\rm d}x}
\def\k{k}

\def\ba{{\bf a}}

\def\bn{{\bf n}}
\def\bp{{\bf p}}

\def\br{{\bf r}}
\def\bu{{\bf u}}
\def\bv{{\bf v}}

\def\bx{{\bf x}}

\def\bB{{\bf B}}

\def\bE{{\bf E}}

\def\F{{\mathsf F}}
\def\G{{\mathsf G}}
\def\I{{\mathsf I}}
\def\K{{\mathbb K}}

\def\M{{\mathsf M}}
\def\N{{\mathsf N}}
\def\P{{\mathsf P}}

\def\R{{\mathsf R}}

\def\W{{\mathsf W}}

\def\cC{{\cal C}}

\def\cE{{\cal E}}
\def\cF{{\cal F}}

\def\cN{{\cal N}}
\def\cP{{\cal P}}

\newcommand{\bref}[1]{(\ref{#1})}

\newtheorem{lemma}{Lemma}[section]
\newtheorem{remark}{Remark}[section]

\def\BeginProof{\underline{Proof.}~}
\def\EndProof{\hfill$\Box$\vspace*{3mm}}

\newcommand{\ASSUM}[2]{(\textbf{#1#2})}

\newcommand{\stabterm}{\textsf{s}}

\newcommand{\bone}{\mathbbm{1}}
\newcommand{\blambda}{{\bm\lambda}}

\begin{document}

\title{{\bf Discretization of mixed formulations of elliptic problems on polyhedral meshes}}

\author{
     Konstantin Lipnikov\footnotemark[1],\quad
     Gianmarco Manzini\footnotemark[1]}

\date{}

\maketitle

\footnotetext[1]{Applied Mathematics and Plasma Physics Group, Theoretical Division,
  Los Alamos National Laboratory, \{lipnikov,manzini\}@lanl.gov}

\begin{abstract}
  We review basic design principles underpinning the construction of
  mimetic finite difference and a few finite volume and finite element
  schemes for mixed formulations of elliptic problems.
  For a class of low-order mixed-hybrid schemes, we show connections
  between these principles and prove that the consistency and
  stability conditions must lead to a member of the mimetic family of
  schemes regardless of the selected discretization framework.
  Finally, we give two examples of using flexibility of the mimetic
  framework: derivation of higher-order schemes and convergent schemes
  for nonlinear problems with small diffusion coefficients.
\end{abstract}

\section{Introduction}

The mixed formulation allows us to calculate simultaneously the
primary solution of a PDE and its flux.
For this reason, mixed formulations are very useful for numerical
solution of multiphysics systems.
The focus of this work is on a single diffusive process that is a part
of almost any complex multiphysics system.

In this paper, we present design principles used in the derivation of
mimetic finite difference (MFD) schemes on polygonal and polyhedral
meshes and establish bridges to design principles used by a few other
discretization frameworks (finite volumes and finite elements).
The focus on the design principle allows us to avoid technical details and 
provide a more clear connection between different frameworks in comparison with
the work performed in \cite{Droniou-Eymard-Gallouet-Herbin:2010}.
We also illustrate the flexibility of the mimetic framework with two
challenging examples: derivation of arbitrary-order accurate schemes
for linear problems and convergent schemes for nonlinear problems with
degenerate diffusion coefficients.

Many ideas underpinning the MFD method were originally formulated in
the sixties for orthogonal meshes using the finite difference
framework from which the name of the method was derived.
Over the years, the MFD method has been extensively developed for the
solution of a wide range of scientific and engineering problems in
continuum mechanics \cite{Margolin-Shashkov-Smolarkiewitcz:2000},
electromagnetics \cite{Hyman-Shashkov:2001,Lipnikov-Manzini-Brezzi-Buffa:2011}, fluid
flows~\cite{Lipnikov-Moulton-Svyatskiy:2008,BeiraodaVeiga-Gyrya-Lipnikov-Manzini:2009,BeiraodaVeiga-Lipnikov-Manzini:2010,Cangiani-Manzini-Russo:2009,BeiraodaVeiga-Droniou-Manzini:2011},
elasticity~\cite{BeiraodaVeiga:2010,BeiraodaVeiga-Mora:2011}, 
obstacle
and control
problems~\cite{AntoniettiBigoniVerani2012,Antonietti-BeiraodaVeiga-Verani:2013,Antonietti-BeiraodaVeiga-Bigoni-Verani:2014-M3AS},
diffusion \cite{Hyman-Shashkov-Steinberg:1997}, discretization of
differential
forms~\cite{Bochev-Hyman:2006,Palha-Rebelo-Hiemstra-Kreeft-Gerritsma:2014-JCP,Brezzi-Buffa-Manzini:2014-JCP},
and eigenvalue analysis~\cite{Cangiani-Gardini-Manzini:2011}.
An extensive list of people who contributed to the development of the
MFD method can be found in the recent book
\cite{BeiraodaVeiga-Lipnikov-Manzini:2014} and review paper
\cite{Lipnikov-Manzini-Shashkov:2014}.
The paper summarizes almost all known results on Cartesian and
curvilinear meshes for various PDEs including the Lagrangian
hydrodynamics.
The book complements the paper by providing numerous examples and
describing basic tools used in the convergence analysis of mimetic
schemes for elliptic PDEs.

The MFD method preserves or mimic essential mathematical and physical
properties of underlying partial differential equations (PDEs) on
general polygonal and polyhedral meshes.
For the elliptic equation, these properties include the \emph{local
  flux balance} and the \emph{duality between gradient and divergence
  operators}.
The latter implies symmetry and positive definiteness of the resulting
matrix operator and is desirable for robustness and reliability of
numerical simulations.
The duality of the primary and derived mimetic operators is one of the
major design principles.
The \emph{definition of the primary mimetic operators is coordinate
  invariant}, which is another design principle that allows us to
build discrete schemes for non-Cartesian coordinate systems.
The discrete operators are also built to satisfy \emph{exact
  identities}, the property that is critical for avoiding spurious
numerical solutions, providing accurate modeling of conservation laws,
and making the convergence analysis possible.

The related discretization frameworks considered in this paper 
include the finite volume methods
\cite{Eymard-Gallouet-Herbin:2010,Droniou-Eymard:2006},
the mixed finite element (MFE) method \cite{Raviart-Thomas:1977}, and the virtual
element method (VEM) \cite{Brezzi:2014}.
Other finite volume frameworks exist that are based on mimetic
principles such as the discrete duality finite volume methods (DDFV),
see, e.g.,~\cite{Domelevo-Omnes:2005,Coudiere-Manzini:2010}, but these
methods do not fit in the MFD framework and will not be considered
here.

The FV methods, originally introduced
in~\cite{Dusinberre:1955,Dusinberre:1961} for the heat equation and
dubbed as the integrated finite difference method, form, perhaps, the
largest class of schemes that can handle unstructured polygonal and
polyhedral meshes, non-linear problems, and problems with anisotropic
coefficients.
An introduction to the finite volume nethoddology can be found in
the recent review \cite{Droniou:2014}.
Almost all FV methods starts with a discrete representation of the
flux balance equation.
This representation is exact and this property is so important that
all the methods that we consider in this paper use the same discrete
form of the balance equation and the difference between them is only
in the discretization of the constitutive equation.

The classical cell-centered FV scheme uses a two-point flux formula
that is second-order accurate for special meshes such as the Voronoi
tessellations.
To overcome this limitation, a class of FV methods, consistent by
design, is proposed by introducing additional unknowns on mesh faces.
Examples of such methods are the \emph{hybrid finite volume} method
\cite{Eymard-Gallouet-Herbin:2010}, and
the \emph{mixed finite volume} method \cite{Droniou-Eymard:2006}.
These FV methods start with different definitions of the cell-based
discrete gradient that are exact for linear solutions.
The formula for the numerical flux based on this gradient needs a
stabilization term.
Construction of the stabilized flux uses two principles.
First, the stabilization term should be zero on linear solutions.
Second, the stabilized flux is defined as the solution of a certain equation with 
a symmetric and positive definite bilinear form.
We will show that these design principles imply the duality principle
in the mimetic framework.

The VEM was originally introduced as an evolution of the MFD method.
In the classical finite element spirit, the duality principle is
incorporated directly in the weak formulation.
The exact identities are replaced by the exact sequence of virtual
finite element spaces.
A new design principle is the unisolvency property where the space of
degrees of freedom is isomorphic to a space of finite element
functions and includes polynomial as well as non-polynomial functions.
The bilinear forms are split explicitly into consistency and stability
forms using problem-dependent $L^2$ and $H^1$ projectors.
We discuss how the new design principles are connected to the
stability and consistency conditions in the mimetic framework.

The recent developments of the MFD framework exploit it flexibility
for selecting non-standard degrees of freedom, optimization of inner
products, and non-standard approximations of primary operators to
build schemes with higher order of accuracy and convergence schemes
for nonlinear PDEs with degenerate coefficients (see also
Section~\ref{sec:divk}).

Extension to higher-order mixed scheme is almost straightforward in
the mimetic framework.
The key step is the proper selection of degrees of freedom that (a)
simplify the discretization of the primary divergence operator and (b)
allows us to formulate a computable consistency condition.
A new design principle is introduced in this case, which states that a
commuting relation exists between the interpolation operators defining
the degrees of freedom of scalar and vector fields, and the divergence
operators in the discrete and continuum settings.

All of the aforementioned methods discretize effectively the
divergence operator ``\emph{div}\,$(\cdot\,)$''.
To solve nonlinear parabolic equations, we employ new MFD schemes
where the primary operator discretizes the combined operator
``\emph{div}\,$\k(\,\cdot\,)$'', where $\k$ is the non-constant
scalar diffusion coefficient.
The resulting scheme uses both cell-centered and face-centered values
of the diffusion coefficient.
This model has applications in heat diffusion \cite{Castor:2004} and
moisture transport in porous media \cite{Richards:1931}.
The duality property mentioned above guarantees that the schemes can
be formulated as algebraic problems with symmetric and positive
definite matrices.
Matrices with these properties lead to better performance of scalable
iterative solvers, such as algebraic multigrid solvers and Krylov
solvers such as the preconditioned conjugate gradient.

Finally, we mention other discretization methods that work on general
meshes.
Our necessarily incomplete list include the polygonal/polyhedral
finite element method
(PFEM)~\cite{%
Wachspress:1975,%
Sukumar-Tabarraei:2004,%
Sukumar-Malsch:2006,%
Tabarraei-Sukumar:2006,%
Manzini-Russo-Sukumar:2014},
hybrid high-order
method~\cite{DiPietro-Ern:2014},
the discontinuous Galerkin (DG) method~\cite{DiPietro-Ern:2011}, 
hybridized discontinuous Galerkin (HDG)
method~\cite{Cockburn-Gopalakrishnan-Lazarov:2009},
and the weak Galerkin (wG) method \cite{Wang-Ye:2014}.

The outline of the paper is as follows.
In Section~\ref{sec:basic:principles:MFD:framework}, we review the
basic discretization principles of the mimetic framework.
In Section~\ref{sec:mixed:mimetic:formulation}, we derive the mimetic
finite difference method for elliptic problems through the consistency and stability properties.
We also prove that any mixed-hydrid method that uses the same degrees
of freedom leads to a member of the mimetic family of schemes.
In Section~\ref{sec:recent:developments}, we review the recent
progress in the development of mimetic methods for mixed formulations
of elliptic problems.
Our final remarks and conclusions are given in
Section~\ref{sec:conclusions}.

\section{Principles of the mimetic discretization framework}
\label{sec:basic:principles:MFD:framework}
\setcounter{equation}{0}

The MFD method mimics important mathematical and physical properties
of underlying PDEs.
We give two examples showing the importance of preserving such
properties in physics simulations.

\begin{figure}[h!]
\centering
\includegraphics[scale=0.75]{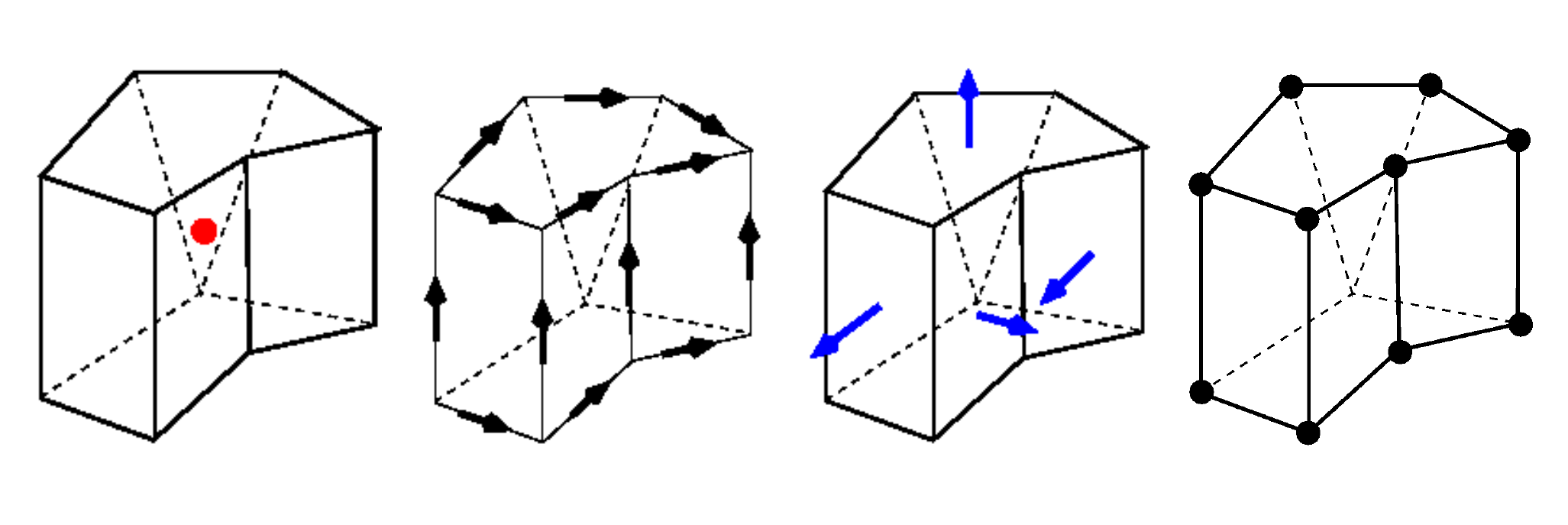}\quad
\caption{Illustration of the degrees of freedom in low-order mimetic
  schemes.  The local spaces associated with cell $c$ from left to
  right are $\cC_{h,c}$, $\cE_{h,c}$, $\cF_{h,c}$, $\cN_{h,c}$. The
  degrees of freedom are shown only on the visible objects.}
\label{fig:dofs}
\end{figure}

Consider a polygonal or polyhedral mesh $\Omega_h$.
We denote the sets of mesh nodes, edges, faces, and cells by symbols
$\cN$, $\cE$, $\cF$ and $\cC$, respectively, the set of vectors
collecting the degrees of freedom associated with those mesh objects
by the corresponding symbol with the subscript ``$h$'', and the
restriction to cell ``$c$'' by the subscript ``$h,c$''.
Each set of vectors of degrees of freedom with the (obvious)
definitions of addition and multiplication by a scalar number is a
linear space.
For example, $\cF_h$ is the linear space of vectors formed by the
degrees of freedom located on the mesh faces, and $\cF_{h,c}$ is its
restriction to cell $c$.
Its precise definition depends on the scheme.
An illustration of particular discrete spaces restricted to a single
cell is shown in Fig.~\ref{fig:dofs}.

The mimetic finite difference method operates with discrete analogs of
the first-order operators.
These operators are designed to satisfy exact identities and duality
principles.

\subsection{Global mimetic operators}

In the mimetic framework, we usually discretize pairs of adjoint
operators, such as the primary divergence $\DIV\colon\cF_h\to\cC_h$
and the derived gradient $\GRADt\colon \cF_h \to \cC_h$.
Hereafter, we will distinguish the derived operators from the primary
operators by using a tilde on the operator's symbol.
It is convenient to think about these operators as matrices acting
between finite dimensional linear spaces.
To discretize a large class of PDEs, we need three pairs of primary
and derived operators, which are discrete analogs of gradient, curl
and divergence operators.
Each pair of operators satisfies a discrete integration by parts
formula, e.g.
\begin{equation}\label{dualGRAD}
  \big[\DIV\, \bu_h,\, q_h\big]_{\cC_h} = -\big[\bu_h,\, \GRADt\, q_h\big]_{\cF_h}
  \qquad \forall \bu_h\in \cF_h,\  \forall q_h \in \cC_h.
\end{equation}
This formula represents one of the mimetic discretization principles
as it mimics the continuum Green formula
\begin{equation}\label{Green}
  \Int_\Omega (\div\, \bu)\, q \dx
  = -\Int_\Omega \bu \cdot \nabla q \dx
  \qquad \forall \bu \in H_{div}(\Omega),\  \forall q \in H^1_0(\Omega).
\end{equation}
The brackets $[\cdot,\cdot]_{\cC_h}$ and $[\cdot,\cdot]_{\cF_h}$
in~\eqref{dualGRAD} stand for an approximation of the integrals
in~\eqref{Green} and will be referred to as the \emph{mimetic inner
products} (or, simply \emph{inner products}).
The inner products are typically constructed from local
(e.g. cell-based or node-based) inner products, which simplifies their
derivation.
For example, the two inner products in~\eqref{dualGRAD} can be
reformulated as
\begin{equation}\label{additivity}
  \big[\bu_h,\, \bv_h\big]_{\cF_h}
  = \Sum\limits_{c \in \Omega_h} \big[\bu_c,\, \bv_c\big]_{\cF_h, c},
  \qquad
  \big[p_h,\, q_h\big]_{\cC_h}
  = \Sum\limits_{c \in \Omega_h} \big[p_c,\, q_c\big]_{\cC_h, c}
\end{equation}
where $\bv_c$, $\bu_c$, $q_c$ and $p_c$ denote the restriction to mesh
cell $c$ of the corresponding global vectors in the left-hand side and
$[\,\cdot,\cdot\,]_{\cF_h, c}$ and $[\,\cdot,\cdot\,]_{\cC_h, c}$ are
the local contribution from $c$ to the global inner products
$[\,\cdot,\cdot\,]_{\cF_h}$ and $[\,\cdot,\cdot\,]_{\cC_h}$,
respectively.

Let $\M_{\cC}$ and $\M_{\cF}$ be the symmetric positive definite
matrices induced by the inner products $[\cdot,\cdot]_{\cC_h}$ and
$[\cdot,\cdot]_{\cF_h}$, respectively.
Then, the explicit formula for the derived gradient operator is
$$
  \GRADt = -\M_{\cF}^{-1}\, \DIV^T\, \M_{\cC}.
$$
This formula shows that this operator has a nonlocal stencil when
matrix $\M_{\cF}$ is irreducible as is typical for unstructured
meshes.
Note that the same property holds true for many other discretization
frameworks.

Formula \eqref{dualGRAD} implies that in general the discrete
operators cannot be discretized independently.
If we discretize one of the operators, e.g., the divergence, and
select the inner products, the other operator, the gradient, must be
derived from formula \eqref{dualGRAD}.  
The existing freedom is in the selection of the inner products.

The selection of the discrete spaces is typically done to simplify the
discretization of the primary mimetic operator.
For a different pair of discrete spaces, e.g., $\cN_h$ and $\cE_h$,
one can find that the gradient operator can be discretized in a
natural way.
In such a case, the gradient operator $\GRAD\colon \cN_h \to \cE_h$ is
the primary mimetic operator and the discrete divergence operator
$\DIVt\colon\cE_h\to\cN_h$ is the derived operator.
This pair of operators satisfies another discrete integration by parts
formula that mimics \eqref{Green}:
\begin{equation}\label{dualDIV}
  \big[\DIVt\, \bv_h,\, p_h\big]_{\cN_h} = -\big[\bv_h,\, \GRAD\, p_h\big]_{\cE_h}
  \qquad \forall p_h\in \cN_h,\  \forall \bv_h \in \cE_h.
\end{equation}
The explicit formula for the derived divergence operator is
$$
  \DIVt = -\M_{\cN}^{-1}\, \GRAD^T\, \M_{\cE},
$$
where $\M_{\cE}$ and $\M_{\cN}$ are symmetric positive definite matrices induced
by the inner products.
Note that in some low-order mimetic schemes, matrix $\M_{\cN}$ is diagonal,
so that the derived operator has a local stencil.

The third pair of discrete operators approximates the continuum
operators ``\emph{curl}''.
Let $\CURL\colon\cE_h\to \cF_h$ and $\CURLt\colon\cF_h\to \cE_h$
satisfy the discrete integration by parts formula
$$
\big[\CURL\, \bv_h,\, \bu_h\big]_{\cF_h} = \big[\bv_h,\, \CURLt\, \bu_h\big]_{\cE_h}
\qquad \forall \bv_h\in \cE_h,\  \forall \bu_h \in \cF_h,
$$
which mimics the continuum formula
$$
\Int_\Omega (\curl\, \bv)\cdot \bu \dx
= \Int_\Omega \bv \cdot (\curl\, \bu) \dx
\qquad \forall \bu \in H_{curl}^0(\Omega),\  \forall \bv \in H_{curl}(\Omega).
$$
The explicit formula for the derived curl operator is
$$
  \CURLt = \M_{\cE}^{-1}\, \CURL^T\, \M_{\cF},
$$
so that this derived operator has typically a non-local stencil.

The spaces of discrete functions that we have introduced so far
satisfy homogeneous boundary conditions.
In \cite{Hyman-Shashkov:1998}, these spaces were enriched conveniently to
approximate the boundary integrals that appear in general Green
formulas.
The resulting derived mimetic operators include an approximation of
the boundary conditions.
We will not follow this approach here, since the focus of this paper
is on mixed-hybrid formulations, which provide another way to
incorporate boundary conditions in a numerical scheme.

The duality of the discrete operators helps us to build numerical schemes
that satisfy discrete conservation laws.
For example, consider the Euler equations in the Lagrangian form:
\begin{equation}
  \label{Euler:Lagrangian}
  \Frac{1}{\rho} \Frac{d \rho}{d t} = -\div\, \bu,
  \qquad
  \rho \Frac{d \bu}{d t} = -\nabla p,
  \qquad
  \rho \Frac{d e}{d t} = -p\,\div\, \bu,
\end{equation}
where $p$ is the pressure, $\rho$ is the density, $\bu$ is the
velocity, and $e$ is the internal energy.
The system is closed by an equation of state.
A mimetic discretization of~\eqref{Euler:Lagrangian} is given by 
\begin{equation}
  \label{Euler:Lagrangian:mimetic}
  \Frac{1}{\rho_h} \Frac{d \rho_h}{d t} = -\DIV\, \bu_h,
  \qquad
  \rho_h \Frac{d \bu_h}{d t} = -\GRADt p_h,
  \qquad
  \rho_h \Frac{d e_h}{d t} = -p_h\,\DIV\, \bu_h,
\end{equation}
where $p_h$, $\rho_h$, $\bu_h$ and $e_h$ are the discrete analogs of
the corresponding continuum quantities that appear
in~\eqref{Euler:Lagrangian} and $\DIV$ and $\GRADt$ are the mimetic
operators acting, respectively, as divergence and gradient.
Let us assume that no external work is done on the system, e.g., $p=0$
of $\partial\Omega$.
The integration by parts and the continuity equation
from~\eqref{Euler:Lagrangian} lead to the conservation of the total
energy $E$: 
\begin{equation}
  \label{energy:conservation}
  \Frac{d E}{d t} 
  = \Int_{\Omega(t)} \rho \Big(\Frac{d \bu}{d t}\cdot \bu +\Frac{d e}{d t}\Big) \dx
  = -\Int_{\Omega(t)} \big(\bu \cdot \nabla p + p\,\div\, \bu) \dx
  = 0.
\end{equation}
To mimic this property, we need the discrete gradient and divergence
operators $\GRADt$ and $\DIV$ to satisfy a discrete integration by
parts formula like~\eqref{dualGRAD}.
Using the same argument that leads to~\eqref{energy:conservation}, we
obtain the conservation of the total discrete energy $E_h$:
\begin{equation*}
  \Frac{d E_h}{d t}
  = -\big[\bu_h,\, \GRADt\, p_h\big]_{\cF_h} - \big[p_h,\,\DIV \bu_h\big]_{\cC_h} = 0.
\end{equation*}
We emphasize that numerical methods that conserve energy usually
have other important properties such as correct prediction 
of a shock position and bounded numerical solution.

Another discretization principle is to derive primary operators that 
mimic exact identities.
This is typically achieved by using the first principles (the divergence
and Stokes theorems) to define the primary operators, e.g. \bref{primary-div}.
As the result, we have
\begin{align*}
  \DIV\, \CURL\, \bv_h = 0,\qquad
  \CURL\, \GRAD\, p_h = 0
  \qquad \forall \bv_h \in \cE_h,\ \forall p_h \in \cN_h.
\end{align*}
Another consequence of the duality principle is that similar
identities hold for the derived mimetic operators.
Using the aforementioned explicit formulas for these operators, we
immediately obtain that
\begin{align*}
  \DIVt\, \CURLt\, \bu_h = 0,\qquad
  \CURLt\, \GRADt\, q_h = 0
  \qquad \forall \bu_h \in \cF_h,\ \forall q_h \in \cC_h.
\end{align*}

These exact identities allows us to design numerical schemes without
non-physical spurious modes.
For instance, in the numerical solution of Maxwell's equations, such
operators guarantee that the magnetic field $\bB_h$ remains
divergence-free for all times.
Applying the primary divergence operator to a semi-discrete form of
Faraday's law of induction, i.e., $\partial\bB_h\slash{\partial
  t}=-\CURL\,\bE_h$, we obtain
\begin{align*}
  \Frac{\partial}{\partial t} \big(\DIV\, \bB_h\big) 
  = \DIV\,\Frac{\partial \bB_h}{\partial t} 
  = -\DIV\,\CURL\,\bE_h
  = 0.
\end{align*}
Therefore, if $\bB_h$ is such that $\DIV\,\bB_h=0$ at time $t=0$, this
relation will be satisfied at any time $t>0$.

\subsection{Local mimetic operators}
\label{subsec:local:mimetic:operators}

From this section, we limit our discussion to elliptic problems and
one pair of the primary and derived operators.
For the practical implementation of mimetic schemes, it is convenient to
write a local integration by parts formula that implies the global
one.
To do it, we need an additional space $\Lambda_h$ of pressure
unknowns defined typically on mesh faces.
Its restriction to cell $c$ is denoted by $\Lambda_{h,c}$ and consists
of the vectors $\blambda_c$.
We recall that the subscript ``$c$'' is added to denote the local
mimetic operators and the local discrete spaces corresponding to cell
$c$.

Let $\DIV_c\colon \cF_{h,c} \to \cC_{h,c}$ be the primary divergence operator.
The derived gradient operator $\GRADt_c\colon \cC_{h,c} \times \Lambda_{h,c}
\to \cF_{h,c}$ satisfies the discrete integration by parts formula
\begin{align}
  &\hspace{2cm}
  \big[\DIV_c\, \bu_c,\, q_c\big]_{\cC_h,c} 
  - [\bu_c,\, \blambda_c]_{\Lambda_{h,c}}
  = -\big[\bu_c,\, \GRADt_c\, 
  \left(\!\!\begin{array}{c}q_c \\ \blambda_c \end{array}\!\!\right)\big]_{\cF_h,c}
  \nonumber\\[0.5em]
  &\hspace{7.5cm}
  \forall \bu_c\in \cF_{h,c},\ 
  \forall q_c\in \cC_{h,c},\ 
  \forall \blambda_c\in\Lambda_{h,c},
  \label{eq:formula:B}
\end{align}
which mimics the continuum Green formula for cell $c$:
$$
  \Int_c (\div\, \bu)\, q \dx
  - \Int_{\partial c} (\bu \cdot \bn)\, q \dx
  = -\Int_c \bu \cdot \nabla q \dx
  \qquad \forall \bu \in H_{div}(c),\  \forall q \in H^1(c).
$$
In order to recover formula \eqref{dualGRAD} for the global discrete
gradient operator, we impose the continuity of $\blambda_c$ and $\bu_c$
on the mesh faces, define the local divergence operator as the
restriction of the global one, define the local spaces as restrictions
of global ones, require that the interface terms cancel each
other,
\begin{equation}\label{assum1}
  \Sum\limits_{c \in \Omega_h} [\bu_c,\, \blambda_c]_{\Lambda_{h,c}} = 0,
\end{equation}
and that the local inner products are summed up into global inner
products as in~\bref{additivity}.

The derivation of the mimetic method follows three generic steps.
First, we select the degrees of freedom such that the local primary
operator, e.g., $\DIV_c$, has a simple form.
Second, we define the inner products in the discrete spaces that
satisfy the consistency and stability conditions.
Third, we postulate the discrete integration by parts formula and
obtain the derived operator, e.g., $\GRADt_c$, from it.
Note that the local derived operator is defined uniquely.

These three steps are discussed in
Section~\ref{sec:mixed:mimetic:formulation} for the mixed formulation
of the diffusion problem.
The flexibility of the mimetic framework is exploited in
Section~\ref{sec:divk}, where we derive another pair of primary
divergence and derived gradient operators for a nonlinear parabolic
problem.
More examples of mimetic schemes can be found in
\cite{BeiraodaVeiga-Lipnikov-Manzini:2014}.

\subsection{Material properties}

The material properties are often included in the definition of the
derived mimetic operator.
Indeed, the Green formula~\eqref{Green} can be rewritten as follows:
\begin{equation}
  \label{Green:material-properties}
  \Int_\Omega (\div\, \bu)\, q \dx
  = -\Int_\Omega \K^{-1} \bu \cdot (\K \nabla) q \dx
  \qquad \forall \bu \in H_{div}(\Omega),\  \forall q \in H^1_0(\Omega).
\end{equation}
According to~\eqref{Green:material-properties}, we can define $\GRADt$
as an approximation of the combined operator $\K\nabla$ and the inner product
$[\bu_h,\bv_h]_{\cF_h}$ as an approximation of the right-hand side
integral $\Int_\Omega \K^{-1} \bu \cdot \bv \dx$, provided that
$\bu_h$ and $\bv_h$ are the degrees of freedom of $\bu$ and $\bv$.

For a perfectly conducting medium, we have the following duality relationship for 
the first-order curl operators:
\begin{equation*}
  \int_{\Omega}\curl\,\bE\cdot\mu^{-1}\bB \dx = 
  \int_{\Omega}\varepsilon\bE\cdot\left(\varepsilon^{-1}\curl\,\mu^{-1}\bB\right)\dx.
\end{equation*}
In this case the inner products in spaces $\cE_h$ and $\cF_h$ are the weighted 
inner products.
The weights are the magnetic permeability $\mu^{-1}$ and electric
permittivity $\varepsilon$.
The derived curl operator $\CURLt$ is an approximation of the combined
operator $\varepsilon^{-1}\curl\,\mu^{-1}$.

\section{Mixed formulation of diffusion problem}
\label{sec:mixed:mimetic:formulation}
\setcounter{equation}{0}

Let $\Omega \in \Re^d$ be a polygonal ($d=2$) or polyhedral ($d=3$) domain 
with the Lipschitz continuous boundary.
Consider the mixed formulation of the diffusion problem:
\begin{equation}
  \begin{array}{rcll}
    \bu         &=& -\K\nabla p & \qquad\textrm{in~}\Omega,\\[1ex]
    \div \, \bu &=& b           & \qquad\textrm{in~}\Omega,
  \end{array}
\end{equation}
subject to the homogeneous Dirichlet boundary conditions on $\partial\Omega$.
As usual, we will refer to the scalar unknown as the \emph{pressure}
and to the vector unknown as the \emph{flux}.
We assume that the diffusion tensor is piecewise constant on mesh
$\Omega_h$ and we denote its restriction to cell $c$ by $\K_c$.
If $\K_c$ is not constant on $c$, we can take its values at the
centroids of the mesh cells without losing the approximation order.

In this section, we consider one local mimetic formulation, two FV
schemes and two mixed-hybrid FE schemes with the same set of degrees
of freedom.
Moreover, we consider the schemes that use the same discrete
divergence operator and the discrete flux balance equation:
\begin{equation}\label{Model}
  \bu_c = {\cal L}_c(p_c, \blambda_c),
  \quad
  \DIV_c \bu_c = b_c^I,
\end{equation}
where ${\cal L}_c$ is a linear operator and $b_c^I$ is defined later.

\subsection{Regular polygonal and polyhedral meshes}\label{sec:mesh}

The analysis of discretization schemes is typically conducted on a
sequence of conformal meshes $\Omega_h$ where $h$ is the diameter of
the largest cell in $\Omega_h$ and $h \to 0$.
A mesh is called conformal if the intersection of any two distinct
cells $c_1$ and $c_2$ is either empty, or a few mesh points, or a few
mesh edges, or a few mesh faces.
Cell $c$ is defined as a closed domain in $\Re^3$ (or $\Re^2)$ with
flat faces and straight edges.

\begin{figure}[ht]
  \begin{center}
    \includegraphics[scale=0.5]{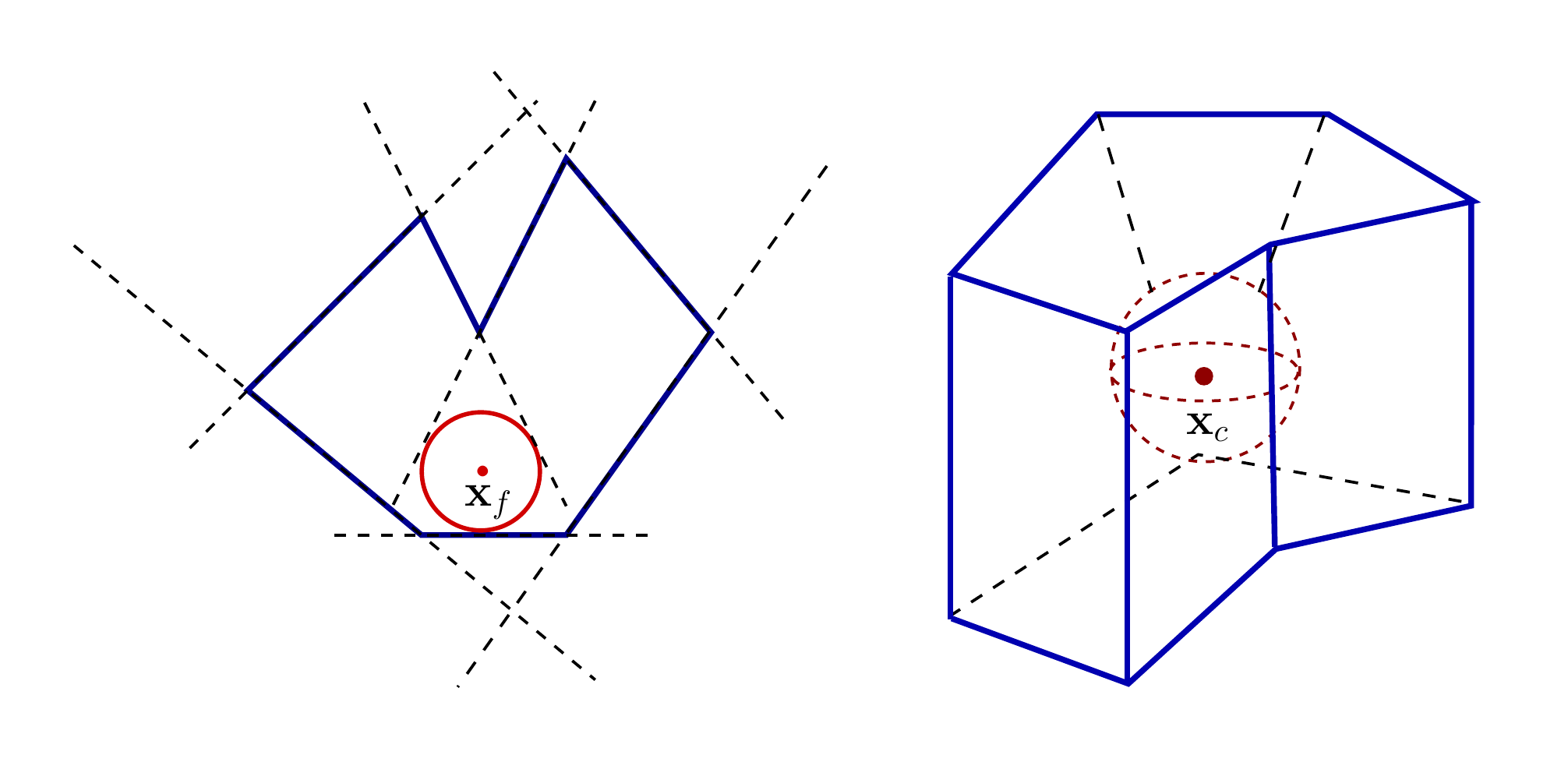}
    \caption{Shape-regular mesh objects.}
    \label{starshaped}
  \end{center}
\end{figure}

Following ~\cite{BeiraodaVeiga-Lipnikov-Manzini:2014}, we make a few
assumptions on the regularity of 3D meshes.
Similar assumptions can be derived for 2D meshes by reducing the
dimension.
Let $n_\star$, $\rho_\star$ and $\gamma_\star$ denote various mesh
independent constants explained below.
\begin{description}
\item\ASSUM{M}{1} 
  Every polyhedral cell $c$ has at most $n_\star$ faces and each face
  $f$ has at most $n_\star$ edges.

\item\ASSUM{M}{2}
  For every cell $c$ with faces $f$ and edges $e$, we have
  \begin{equation}\label{shape}
    \rho_\star\, \big({\rm diam}(c)\big)^3 \le\,|c|, \quad\,
    \rho_\star\, \big({\rm diam}(c)\big)^2 \le\,|f|, \quad\,
    \rho_\star\, {\rm diam}(c) \le\,|e|,
  \end{equation}
  where $|\cdot|$ denotes the Euclidean measure of a mesh object.

\item\ASSUM{M}{3}
  For each cell $c$, there exists a point $\bx_c$ such that $c$ is
  star-shaped with respect to every point in the sphere of radius
  $\gamma_\star\,{\rm diam}(c)$ centered at $\bx_c$.
  For each face $f$, there exists a point $\bx_f\in f$ such that $f$
  is star-shaped with respect to every point in the disk of radius
  $\gamma_\star\, {\rm diam}(c)$ centered at $\bx_f$ as shown in
  Fig.~\ref{starshaped}.

\item\ASSUM{M}{4} 
  For every cell $c$, and for every $f\in \partial c$, there exists a
  pyramid contained in $c$ such that its base equals $f$, its
  height equals ${\gamma_\star}\,{\rm diam}(c)$ and the projection
  of its vertex onto $f$ is $\bx_f$.
\end{description}

The conditions \ASSUM{M}{1}-\ASSUM{M}{4} are sufficient to develop an
{\it a priori} error analysis of various discretization schemes.  We
recall only two results underpinning this error analysis.
The first one is the Agmon inequality that uses \ASSUM{M}{4} and
allows us to bound traces of functions:
\begin{equation}
\label{eq:Agmon}
  \sum_{f\in\partial c}\|q\|_{L^2(f)}^2
  \leq C \left(
     \big({\rm diam}(c)\big)^{-1}\|q\|_{L^2(c)}^2 + {\rm diam}(c)\, |q|_{H^1(c)}^2
  \right)
  \qquad \forall q \in H^1(c).
\end{equation}
The second one is the following approximation result: For any function
$q \in H^2(c)$ there exists a polynomial $q^1 \in \cP^1(c)$ such that
\begin{equation}
  \|q - q^1\|_{L^2(c)} +
  {\rm diam}(c)\, |q-q^1|_{H^1(c)}
  \leq
  C \big({\rm diam}(c)\big)^2|q|_{H^2(c)}.
\label{eq:approximation}
\end{equation}
Hereafter, we will use symbols $C$, $C_1$, $C_2$ to denote generic
constants independent of $h$.

\begin{remark}
  The mesh regularity assumptions \ASSUM{M}{1}--\ASSUM{M}{4} are not unique.
  They could be generalized to non-star shaped cells
  \cite{BeiraodaVeiga-Lipnikov-Manzini:2014} by splitting each
  polygonal or polyhedral cell into the finite number of shape regular simplexes.
\end{remark}

\subsection{Mimetic discretization framework}
\label{subsec:mimetic:discretization:framework}

As pointed out at the end of
subsection~\ref{subsec:local:mimetic:operators}, the \emph{first step}
of the construction of a mimetic scheme consists in the selection of
the degrees of freedom.
The degrees of freedom are such that the primary divergence operator
has a simple form.

The discrete space $\cC_h$ consists of one degree of freedom per cell;
its dimension equals the number of mesh cells; and for each vector
$p_h\in\cC_h$ we shall denote the value of $p_h$ associated with cell
$c$ by $p_c\in\cC_{h,c}$
Furthermore, we denote the vector of degrees of freedom of a smooth
function $p$ by $p_h^I \in \cC_h$.
In the MFD method $p_c^I$ is typically the cell average of $p$ over cell $c$
and this definition is used by the existing superconvergence
analysis~\cite{Brezzi-Lipnikov-Shashkov:2005}.
%
%

The discrete space $\Lambda_h$ consists of one degree of freedom per
mesh face, e.g., $\lambda_f$; its dimension equals the number of mesh
faces; and for each vector $\lambda_h\in\Lambda_h$ we shall denote its
restriction to cell $c$ by $\blambda_c\in\Lambda_{h,c}$.
The continuity of local vectors $\blambda_c$ across mesh cells is
satisfied automatically.
The value $\lambda_f$ can be associated with the value of a smooth
function $p$ at the face centroid.

The discrete space $\cF_h$ consists of one degree of freedom per
boundary face and two degrees of freedom per interior face.
For vector $\bu_h \in \cF_h$, we denote by $\bu_c$ its restriction to
cell $c$, and by $u^c_f$ its component associated with face $f$ of
cell $c$.
For a smooth function $\bu$, we denote by $\bu_h^I \in \cF_h$ the
vector of degrees of freedom.
The value $(u^c_f)^I$ is defined as the integral average of flux $\bu
\cdot \bn_f$ through face $f$, where $\bn_f$ is the face normal fixed once
and for all.
Hereafter, we consider a subspace of $\cF_h$ whose members satisfy the
flux continuity constraint
\begin{equation}\label{continuity}
  u^{c_1}_f = u^{c_2}_f
\end{equation}
on each interior face $f$ shared by cells $c_1$ and $c_2$.
With a slight abuse of notation, we shall refer to $\cF_h$ as the
space that satisfies condition \eqref{continuity}.

The local primary divergence operator is defined using a
straightforward discretization of the divergence theorem:
\begin{equation}\label{primary-div}
  \left.(\DIV \bu_h)\right|_c 
  \equiv\DIV_c\bu_c 
  = \frac{1}{|c|} \Sum\limits_{f\in\partial c}|f|\,\sigma_{c,f}\, u^c_f,
\end{equation}
where $\sigma_{c,f}$ is either $1$ or $-1$ depending on the mutual
orientation of the fixed normal $\bn_f$ and the exterior normal
$\bn_{c,f}$ to $\partial c$.
Observe that this definition remains the same in all coordinate systems.

\medskip 
The \emph{second step} of the construction of a mimetic scheme is to
define accurate inner products in $\cC_{h,c}$ and $\cF_{h,c}$ that
satisfy the consistency and stability conditions.
The inner product for space $\cC_{h,c}$ is simple:
\begin{align}
  \big[p_c,\, q_c\big]_{\cC_{h,c}} = |c|\, p_c\, q_c. 
  \label{eq:inner:product:scalar:low:order}
\end{align}
Let ${\cal SC}_c$ be the space of constant functions.  
Then, the above inner product implies the obvious result:
$$
  \big[p_c^I,\, q_c^I\big]_{\cC_{h,c}} = \int_c p\, q \dx
  \qquad\forall p\in \cP^0(c),\ \forall q \in {\cal SC}_c.
$$
Despite its simplicity, we can use this relation to formulate the
general principle that we apply to the derivation of other inner
products.
We define the consistency condition as the following exactness
property: 
\emph{the $L^2$ inner product of two smooth functions is equal to the
  mimetic inner product of their interpolants when one of the
  functions ($p$ in this case) is a polynomial of a given degree and
  the other one belongs to a sufficiently rich space (possibly
  infinite dimensional) that must include polynomial functions.}

The exactness property implies that $L^2$ inner products of a large
class of functions can be calculated exactly using the degrees of
freedom.
In particular, we write the consistency condition for
the inner product on space $\cF_{h,c}$ as the exactness property:
\begin{equation}\label{consistency}
  \big[\bu_c^I,\, \bv_c^I\big]_{\cF_{h,c}} = \int_c \K_c^{-1} \bu\cdot \bv \dx
  \qquad\forall \bu\in (\cP^0(c))^{d},\ \forall \bv \in {\cal SF}_c,
\end{equation}
where ${\cal SF}_c$ is a specially designed space containing the constant vector-functions:
$$
{\cal SF}_c = \left\{ \bv\colon\ 
  \div\,\bv \in \cP^0(c),\quad 
  \bv\cdot \bn_f \in \cP^0(f)\quad \forall f \in \partial c \right\}
$$
Let us show that the right-hand side of \bref{consistency} can be
calculated using the degrees of freedom.
Let $q^1$ be a linear function such that $\K_c\,\nabla q^1 = \bu$.
Inserting $\bu$ in~\bref{consistency}, integrating by parts, and
using the properties of ${\cal SF}_c$, we obtain
\begin{align*}
  \Int_c \nabla q^1\cdot \bv \dx
  = -\Int_c (\div\,\bv) q^1 \dx + \Int_{\partial c} (\bv \cdot \bn)\, q^1 \dx
  = \sum\limits_{f\in \partial c}\bv\cdot\bn_{c,f}\left(
    \Int_f q^1 \dx - \frac{|f|}{|c|}\Int_c q^1 \dx
  \right).
\end{align*}
Let $\M_{\cF,c}$ be the inner product matrix and
$\br_c(q^1)\in\cF_{h,c}$ be the vector with components 
$\sigma_{c,f}\left(\Int_f q^1\dx-\frac{|f|}{|c|}\Int_c q^1\dx\right)$.
Then, combining the last formulas, we have
\begin{align*}
  \big[(\K_c\,\nabla q^1)_c^I,\, \bv_c^I\big]_{\cF_{h,c}} 
  &= \big((\K_c\,\nabla q^1)_c^I\big)^T\, \M_{\cF,c}\, \bv_c^I 
  = \sum\limits_{f\in \partial c} \bv\cdot\bn_f\sigma_{cf}
  \left(
    \Int_f q^1\dx - \frac{|f|}{|c|}\Int_c q^1\dx
  \right)
  \\
  &= \big(\br_c(q^1)\big)^T \bv_c^I.
\end{align*}
Since $\bv$ is any function in ${\cal SF}_c$, we can show that its
interpolant $\bv_c^I$ is any vector in $\cF_{h,c}$.
Indeed, it is sufficient to define a few functions $\bv$ as the solutions 
of cell-based PDEs with different boundary conditions.
Hence, the consistency condition gives us the following matrix
equations
\begin{equation}
  \label{matrix-equations}
  \M_{\cF,c}\, (\K_c\,\nabla q^1)_c^I = \br_c(q^1)
  \qquad\forall q^1 \in \cP^1(c).
\end{equation}
Due to the linearity of these equations, it is sufficient to consider
only three (two in two dimensions) linearly independent linear
functions: $q^1_x = x$, $q^1_y = y$, and $q^1_z = z$.
Let
\begin{align}
  \N_c = [(\K_c\,\nabla x)_c^I\ (\K_c\,\nabla y)_c^I\ (\K_c\,\nabla z)_c^I],
  \quad 
  \R_c = [\br_c(x)\ \br_c(y)\ \br_c(z)]
\label{eq:N:R:general:definition}
\end{align}
be two rectangular $n_c \times 3$ matrices
where $n_c$ is the number of faces in cell $c$.
The equations \bref{matrix-equations} are now equivalent to the matrix 
equation:
\begin{equation}
\label{MN=R}
  \M_{\cF,c} \N_c = \R_c.
\end{equation}
Note that a symmetric positive definite solution $\M_{\cF,c}$ (if it exists) 
is not unique even for a tetrahedral cell $c$.
The existence of solutions is based on the following result proved in 
\cite{BeiraodaVeiga-Lipnikov-Manzini:2014}.
\begin{lemma}\label{lem:NtR}
  Let matrices $\N_c$ and $\R_c$ be defined as in~\bref{eq:N:R:general:definition}.
  Then,
  $
  \R_c^T \N_c = |c|\,\K_c.
  $
\end{lemma}

This lemma allows us to write the explicit formula for matrix
$\M_{\cF,c}$:
$$
  \M_{\cF,c} = \R_c (\R_c^T \N_c)^{-1} \R_c^T + \gamma_c \P_c,
  \quad \P_c = \I - \N_c (\N_c^T \N_c)^{-1} \N_c^T
$$
with a positive factor $\gamma_c$ in front of the projection matrix
$\P_c$.  
A recommended choice for $\gamma_c$ is the mean trace of the first
term.
A family of mimetic schemes is obtained if we replace $\gamma_c$ by an
arbitrarily symmetric positive definite matrix $\G_c$:
\begin{equation}\label{M-family}
  \M_{\cF,c} = \R_c (\R_c^T \N_c)^{-1} \R_c^T + \P_c\, \G_c\, \P_c.
\end{equation}
The stability of the resulting mimetic method depends on the spectral
bounds of matrix $\G_c$ which should be uniformly bounded by $\gamma_c$.
Algebraically, it means that there exists two generic constants $C_1$ 
and $C_2$ such that
\begin{align}
  C_1 |c|\sum_{f\in\partial c}|v^c_f|^2
  \le \bv_c^T\,\M_{\cF,c}\,\bv^T_c \le
  C_2 |c|\sum_{f\in\partial c}|v^c_f|^2
  \label{eq:stab:cond}
\end{align}
holds for every $\bv_c=\{v^c_f\}_{f\in\partial c}$.
This formula is called the \emph{stability condition} in the mimetic discretization framework.
The existence of the mesh independent constants $C_1$ and $C_2$ can be
shown using assumptions \ASSUM{M}{1}-\ASSUM{M}{4}.

\medskip
The \emph{third} step of the construction of a mimetic scheme is
to postulate either the local or global integration by parts formula
and derive the gradient operator from it. 
For instance, given the global discrete operators
$\DIV\colon\cF_h\to\cC_h$ and $\GRADt\colon\cC_h\to\cF_h$, we can
write the mimetic scheme as follows: Find $\bu_h \in \cF_h$ and $p_h \in \cC_h$
such that
$$
  \bu_h = -\GRADt\, p_h,
  \qquad
  \DIV\, \bu_h = b^I,
$$
where $b^I \in C_h$.

The local mimetic formulation requires us to define $[\bu_c,\, \blambda_c]_{\Lambda_{h,c}}$
which satisfies condition \bref{assum1}.
Let
\begin{equation}\label{interface}
  [\bu_c,\, \blambda_c]_{\Lambda_{h,c}} 
  = \sum_{f\in\partial c}\sigma_{c,f}\,|f|\,\lambda_f\,u^c_f.
\end{equation}
Then, the local formulation is to find $\bu_c \in \cF_{h,c}$ and $p_c \in \cC_{h,c}$
in all mesh cells such that
$$
  \bu_c = -\GRADt_c \left(\!\!\begin{array}{c} p_c \\ \blambda_c \end{array}\!\!\right),
  \qquad
  \DIV_c\, \bu_c = b_c^I
$$
subject to the flux continuity conditions \bref{continuity} and the homogeneous 
Dirichlet boundary conditions $\lambda_f = 0$ for $f \in \partial \Omega$.
The local derived operator has the explicit form:
\begin{equation}\label{gradient-cell}
  \GRADt_c \left(\!\!\begin{array}{c}p_c \\ \blambda_c\end{array}\!\!\right)
  = -\M_{\cF,c}^{-1} \left(\!\begin{array}{c} 
    \sigma_{c,f_1}\, |f_1|\,(p_c - \lambda_{f_1}) \\ 
    \vdots \\
    \sigma_{c,f_{n_c}}\,|f_{n_c}|\,(p_c - \lambda_{f_{n_c}})  \end{array}\!\right).
\end{equation}

\begin{lemma}
  Under assumption \bref{additivity} the local and global mimetic
  formulations are equivalent.
\end{lemma}
\BeginProof
Let $\bv_h$ be an arbitrary vector in $\cF_h$ and $\bv_c$ be its
restriction to cell $c$.
To show that the solution to the local mimetic formulation is the
solution to the global one, we first multiply both sides of the local
constitutive equations by $\bv_c$, then sum up the results over the
mesh cells, cancel all internal face terms containing $\lambda_{f}$,
and finally use the local duality relation~\bref{eq:formula:B} between
$\DIV_c$ and $\GRADt_c$ to obtain:
\begin{equation*}
  \sum\limits_c \big[\bu_c,\, \bv_c\big]_{\cF_h,c}
  = -\sum\limits_c 
    \big[\GRADt_c \left(\!\!\begin{array}{c} p_c \\ \blambda_c \end{array}\!\!\right),\, \bv_c\big]_{\cF_h,c}
  = \sum\limits_c \big[\DIV_c\,\bv_c,\, p_c\big]_{\cC_h,c}.
\end{equation*}
From the additivity of the inner products and~\bref{dualGRAD} we
obtain that
\begin{align*}
  \big[\bu_h,\, \bv_h\big]_{\cF_h}
  = \big[\DIV\,\bv_h,\, p_h\big]_{\cC_h}
  = -\big[\GRADt\, p_h,\, \bv_h\big]_{\cF_h}
  \qquad\forall\bv_h\in\cF_h,
\end{align*}
which implies that $\bu_h = -\GRADt\, p_h$.

To show the opposite statement, we repeat the above argument in the
reverse order.
This proves the assertion of the lemma.
\EndProof

Let us consider the matrix equation $\N_c=\W_{\cF,c}\,\R_c$ (compare with~\bref{MN=R}).
To implement the mimetic scheme in a computer program, we need to know only matrix $ \W_{\cF,c}$.
The general solution to the matrix equation is
$$
  \W_{\cF,c} 
  = \N_c\, (\N_c^T \R_c)^{-1}\, \N_c^T 
  + \widetilde\G_c\, (\I - \R_c\,(\R_c^T\R_c)^{-1}\,\R_c^T),
$$
where $\widetilde \G_c$ is an arbitrary $n_c\times n_c$ matrix, possibly non-symmetric.
This formula leads to a large family of stable and unstable schemes that we
refer to as the \emph{extended mixed-hybrid family} of schemes.
Positive definiteness of $\W_{\cF,c}$ is necessary for proving the
scheme's convergence following the path described in
\cite{Lipnikov-Shashkov-Yotov:2009}.

If in addition $\W_{\cF,c}$ is symmetric, then it is one of the
matrices $\M_{\cF,c}^{-1}$.
The classical two-point flux FV scheme is obtained when all matrices
$\W_{\cF,c}$ are positive definite and diagonal.

\subsubsection{Error estimates}

Let $\Omega$ have a Lipschitz continuous boundary.
Furthermore, let every cell $c$ be shape regular as explained in Sec.~\ref{sec:mesh}. 
We assume that $\bx_c$ is the centroid of cell $c$.
We use the triple-bar notation, e.g., $|||\cdot|||$, for the norms
induced by the mimetic inner products.
Then, the interpolants of the exact solution, $p^I \in \cC_h$ and
$\bu^I \in \cF_h$, satisfy \cite{Brezzi-Lipnikov-Shashkov:2005}
$$
  |||p^I - p_h |||_{\cC_h} + ||| \bu^I - \bu_h |||_{\cF_h} \le C\, h,
$$
where $p_h$ and $\bu_h$ are solutions of the global mimetic
formulation.
If in addition $\Omega$ is convex and $C_1$ in \bref{eq:stab:cond} is 
sufficiently large, then
$$
  ||| p^I - p_h |||_{\cC_h} \le C\, h^2.
$$

\subsection{Finite volume discretization framework}

Two examples of FV schemes that fit within the mimetic framework are
the mixed finite volume (MFV) method~\cite{Droniou-Eymard:2006} and
the hybrid finite volume (HFV)
method~\cite{Eymard-Gallouet-Herbin:2010}.
Both methods give a family of schemes because a stabilization
term, which can be suitably parameterized, appears in their
formulation.
Their design principles are based on conditions that imply the 
mimetic duality principle.

Both methods use the same degrees of freedom for pressure and flux
as the local mimetic formulation and the same discrete flux balance
equation, i.e., $\DIV_c\,\bu_c=b_c^I$.
Recall that $\bu_c$ collects the flux unknowns associated with cell
$c$.  
The local pressure unknowns are the cell pressure $p_c$, which is
associated with the interior of $c$, and the interface pressure
$\lambda_f$, which is associated with face $f$.
In the formulation of these FV methods, $p_c$ can be associated with
the pressure value at any point inside the cell.
As in the previous subsection, $\blambda_c=\{\lambda_f\}_{f\in\partial
  c}$ is the vector whose size is equal to the number of faces in cell
$c$.

In the HFV method, a discrete gradient in cell $c$ is defined by
applying the mid-point quadrature rule to the divergence theorem:
\begin{equation}\label{grad-hfv}
  \nabla_c \left(\!\!\begin{array}{c} p_c \\ \blambda_c\end{array}\!\!\right) 
   = \frac{1}{|c|}\sum_{f\in\partial c}|f|\,\lambda_f\,\bn_{c,f}
   = \frac{1}{|c|}\sum_{f\in\partial c}|f|(\lambda_f - p_c)\bn_{c,f}.
\end{equation}
This formula provides the exact value of the gradient whenever $p$ is
a linear function since for any constant vector $\ba$ and any position
vector $\bx_c$, we have the geometric identity
\begin{align*}
  |c|\,\ba = \sum_{f\in\partial c}\,|f|\,\ba\cdot(\bx_f-\bx_c)\,\bn_{c,f},
\end{align*}
where $\bx_f$ is the face centroid.

Formula \bref{grad-hfv} is used to define the numerical scheme after
an additional \emph{stabilization term} $\stabterm_{c,f},$ is
included in the definition of the numerical flux $u_f^c$:
\begin{equation}\label{flux-hfv}
  u^{c}_{f} 
  = -\bn_f\,\cdot \K_c \nabla_c \left(\!\!\begin{array}{c} p_c \\ \blambda_c\end{array}\!\!\right)
    + \stabterm_{c,f}.
\end{equation}
Like in the mimetic framework, the stabilization term is designed very carefully 
to preserve the consistency of the scheme.
Let $\F_c$ be the diagonal matrix with entries $|f|$ on the main
diagonal, $f \in \partial c$.
Similarly, let $\Sigma_c$ be the diagonal matrix with entries $\sigma_{c,f}$.
Furthermore, let $\bone^T=(1,1,\ldots,1)^T$, and
\begin{align*}
  S_{c,f}(p_c,\, \blambda_c) 
  = \lambda_f - p_c - \nabla_c 
    \left(\!\!\begin{array}{c} p_c \\ \blambda_c\end{array}\!\!\right)  \cdot (\bx_f - \bx_c).
\end{align*}
Note that the last expression is zero on linear pressure functions.
The numerical flux is defined implicitly as the solution of
\begin{align}
  \label{eq:FV:stab:term:def}
  (\tilde p_c \bone - \tilde \blambda_c)^T\, \Sigma_c\,\F_c\, \bu_c
  = |c|\, \K_c \nabla_c \left(\!\!\begin{array}{c} p_c \\ \blambda_c\end{array}\!\!\right)
  \cdot \nabla_c \left(\!\!\begin{array}{c} \tilde p_c \\ \tilde\blambda_c\end{array}\!\!\right)
  + \Sum_{f \in \partial c} \alpha_{c,f} \Frac{|f|}{d_{c,f}}\, S_{c,f}(p_c,\blambda_c)\, S_{c,f}(\tilde p_c,\tilde\blambda_c)
\end{align}
for any $\tilde p_c$ and $\tilde\blambda_c$.
Here $\alpha_{c,f}$ is a positive parameter and $d_{c,f}$ is the
distance between $\bx_c$ and the plane containing face $f$.
After a few algebraic manipulations, we obtain:
\begin{align*}
  \stabterm_{c,f}
  = \sigma_{c,f}\sum_{f' \in \partial c}\alpha_{cf'} \frac{|f'|}{d_{cf'}} S_{cf'}(p_c,\blambda_c)
  \left[-\Frac{\delta_{f,f'}}{|f|} + \frac{1}{|c|}\bn_{c,f}\cdot(\bx_{f'}-\bx_{c})\right]
\end{align*}
where $\delta_{f,f'}$ is the Kronecker symbol.
It is obvious that the stabilization term is zero when it is calculated using 
the degrees of freedom of a linear function.

The right-hand side of \bref{eq:FV:stab:term:def} is a symmetric
bilinear form with respect to the pressure unknowns.
It is positive definite when $\alpha_{c,f}> 0$ and $\lambda_f = 0$ on $f \in \partial \Omega$.
It is uniformly
bounded from below when $\alpha_{c,f}$ are sufficiently large and the
mesh satisfies the regularity conditions described in
Section~\ref{sec:mesh}.

\begin{remark}
  Equation~\bref{eq:FV:stab:term:def} is the key design principle. 
  When $p_c$ is associated with the cell centroid, 
  the left-hand side of this equation coincides with
  the left-hand side of~\bref{eq:formula:B} where the divergence is
  given by \bref{primary-div} and the interface term is defined by \bref{interface}.
  So, equation~\bref{eq:FV:stab:term:def} is also a representation of the duality
  principle.
\end{remark}

Formula \bref{grad-hfv} can be rewritten using the mimetic
matrix $\N_c$ as follows:
\begin{align*}
  \nabla^{HFV}_c \left(\!\!\begin{array}{c} p_c \\ \blambda_c\end{array}\!\!\right)
  = \frac{\K_c^{-1}}{|c|}\mathsf{N}_c^T\, \F_c\, \blambda_c.
\end{align*}
This formula should be compared with the formula for the discrete
gradient in the MFV method:
\begin{align*}
  \nabla_c^{MFV} \left(\!\!\begin{array}{c} p_c \\ \blambda_c\end{array}\!\!\right) 
  = -\frac{\K_c^{-1}}{|c|}\,\R_c^T\,\Sigma_c\,\bu_c.
\end{align*}
This formula is exact for a linear pressure function and constant flux,
but it also has to be stabilized.
The vector of numerical fluxes $\bu_c$ is defined as the solution of
\begin{align}
  (\bv_c)^T \G_c \bu_c = (p_c\bone - \blambda_c)^T\, \Sigma_c\, \F_c\, \bv_c
  \qquad \forall \bv_c \in \cF_{h,c}
  \label{eq:flux:def}
\end{align}
where $\G_c$ is a symmetric positive definite matrix.
This design principle shows even clear connection with the mimetic duality principle
\bref{eq:formula:B}.
Since $\mathsf{G}_c$ is invertible, we have that $\bu_c$ is a linear combination of 
$p_c$ and $\blambda_c$, i.e. condition \bref{Model}. 
We can summarize the above discussions in the following lemma.

\begin{lemma}\label{lem:equivalence}
  Any linearity preserving mixed-hybrid scheme of type \eqref{Model}
  with the degrees of freedom given by $\cC_h$, $\Lambda_h$, and $\cF_h$ is a member
  of the extended mixed-hybrid family of schemes.
  In addition, let $\bu_c = {\cal L}_c(p_c,\blambda_c)$ and assume that
  the bilinear form
  \begin{equation}\label{spd}
    {\cal B}\big( (\tilde{p}_c,\tilde{\blambda}_c), (p_c,\blambda_c) \big) :=
    (\tilde p_c \bone - \tilde\blambda_c)^T\,\Sigma_c\,\F_c\,{\cal L}(p_c,\blambda_c)
  \end{equation}
  is symmetric, uniformly coercive and uniformly bounded (with respect
  to some norm that may depend on the problem).
  Then, the resulting scheme belongs to the mimetic family of schemes.
\end{lemma}

\BeginProof
The most general form of the constitutive equation in \bref{Model} is given by the
linear relationship
$$
  \bu_c 
  = {\cal L}_c(p_c,\blambda_c) 
  = -\widetilde \W_{\cF,c} \left(\!\begin{array}{c} \blambda_c \\ p_c\end{array}\!\right)
$$
with a rectangular $n_c \times (n_c + 1)$ matrix $\widetilde \W_{\cF,c} = [\widetilde
\W_{\cF,c}^{(1)},\, \widetilde \W_{\cF,c}^{(2)}]$.
Since the scheme is exact for constant pressure functions, we have
$$
  0 = -\widetilde \W_{\cF,c} \left(\!\begin{array}{c} \bone \\ 1 \end{array}\!\right).
$$
Multiplying the last equation by $p_c$ and subtracting from the previous
one, we obtain
$$
  \bu_c = \widetilde \W_{\cF,c}^{(1)} (p_c\bone - \blambda_c).
$$
By our assumption, this formula is exact for all linear pressure functions
and the corresponding constant flux functions.
Let us take linearly independent pressure functions $x$, $y$, $z$
and use the matrix notations introduced above to derive the following 
consequence of the linearity preservation property:
$$
  \N_c = \widetilde \W_{\cF,c}^{(1)}\, \Sigma_c\,\F_c^{-1}\, \R_c.
$$
Hence, $\widetilde \W_{\cF,c}^{(1)} = \W_{\cF,c}\, \F_c\, \Sigma_c$ and 
the resulting scheme belongs to the extended mixed-hybrid family of schemes.
Now
$$
  (\tilde p_c \bone - \tilde\blambda_c)^T\,\Sigma_c\,\F_c\,\bu_c
  = (\tilde p_c \bone-\tilde\blambda_c)^T\,\Sigma_c\,\F_c\, \W_{\cF,c}\, \F_c\, \Sigma_c (p_c\bone - \blambda_c).
$$
Our symmetry and coercivity assumptions imply that $\W_{\cF,c}$ is
symmetric and positive definite; hence, it is invertible and
$\W_{\cF,c}^{-1}$ coincides with one of the mimetic inner product
matrices $\M_{\cF,c}$.
Let $\tilde \bu_c = \W_{\cF,c}^{(1)} (\tilde p_c\bone - \tilde\blambda_c)$. Then,
$$
 (\tilde p_c \bone - \tilde\blambda_c)^T\,\Sigma_c\,\F_c\,\bu_c
 = \tilde \bu_c^T\, \W_{\cF,c}^{-1}\, \bu_c.
$$
The uniform coercivity and boundness conditions imply that matrix
$\W_{\cF,c}^{-1}$ satisfies the mimetic stability condition.
This proves the assertion of the lemma.
\EndProof

A detailed comparison of MFD, HFV and MFV methods is performed in
\cite{Droniou-Eymard-Gallouet-Herbin:2010}.
In particular, the authors have shown that all schemes can be
generalized to provide the identical families of numerical schemes.

\subsection{Finite element discretization framework}

\subsubsection{Raviart-Thomas FE method on simplexes}
Let us consider the family of mimetic schemes, where the mass matrix
is given in the form of equation~\eqref{M-family}.
When the mesh is formed by simplexes, e.g., triangles in 2D and
tetrahedra in 3D, there is a choice of the stabilization matrix
$\P_c\, \G_c \,\P_c$ that provides the mass matrix from the lowest
order Raviart-Thomas space \cite{Cangiani-Manzini:2008}.
For a simplex, the stabilization matrix is a one-rank matrix, $\P_c\,
\G_c\,\P_c=\bp_c^T\,\big[g_c\big]\,\bp_c$, where
\begin{align*}
  g_c = \frac{1}{d^2(d+1)}\sum_{f\in\partial c}(\bx_f-\bx_c)^T\K_c^{-1}(\bx_f-\bx_c)
\end{align*}
and vector $\bp_c^T = \big(|f_1|, |f_2|, \ldots, |f_{d+1}|\big)$, $f_i \in \partial c$.

\subsubsection{Virtual element method}

The VEM was originally introduced as an evolution of the MFD method.
In the classical finite element spirit, the duality principle is now
incorporated in the weak formulation:
\emph{Find $\bu_h \in \widetilde{\cal
  SF}_h$ and $p_h \in \widetilde \cC_h$ such that}
\begin{subequations}
  \begin{align}
    a_h(\bu_h,\, \bv_h)  - (\div_h\,\bv_h, p_h)_{L^2(\Omega)} &=0\phantom{(b, q_h)_{L^2(\Omega)}} \quad\forall \bv_h \in \widetilde{\cal SF}_h,\label{eq:VEM:1}\\[1ex]
    (\div_h\,\bu_h, q_h)_{L^2(\Omega)}                        &=(b, q_h)_{L^2(\Omega)}\phantom{0} \quad\forall q_h \in \widetilde \cC_h,\label{eq:VEM:2}
  \end{align}
\end{subequations}
where $\widetilde \cC_h$ is the space of piecewise constant functions
isometric to $\cC_h$, $\widetilde{\cal SF}_h$ is a virtual space built
from the local virtual spaces $\widetilde{\cal SF}_c$,
$$
  a_h(\bu_h,\, \bv_h) = \Sum\limits_{c \in \Omega_h} a_{h,c}(\bu_h,\, \bv_h),
  \quad
  a_{h,c}(\bu_h,\, \bv_h) = a_{h,c}^{(1)}(\bu_h,\, \bv_h) + a_{h,c}^{(2)}(\bu_h,\, \bv_h),
$$
$a_{h,c}^{(1)}$ and $a_{h,c}^{(2)}$ are the consistency and stability
terms (to be defined later in the section), and operator $\div_h$ is
defined cell-by-cell as the local $L^2$-orthogonal projection of the
continuum divergence operator onto $\widetilde \cC_h$.
The virtual space $\widetilde{\cal SF}_c$ includes polynomial as well
as non-polynomial functions.

A new design principle is given by the unisolvency property.
The virtual space $\widetilde{\cal SF}_c$ is build to be isomorphic
to the mimetic space $\cF_{h,c}$, i.e., the space of degrees of freedom.
The virtual space is defined as the subspace of ${\cal SF}_c$:
$$
\widetilde{\cal SF}_c = \left\{ \bv\colon\ 
  \div\,\bv \in \cP^0(c),\quad 
  \bv\cdot \bn_f \in \cP^0(f)\quad \forall f \in \partial c,\quad
  \curl\,\bv = 0 \right\}.
$$

In order to describe the splitting of the bilinear form
$a_{h,c}(\bu_h,\, \bv_h)$ we need to introduce the problem-dependent
$L^2$ projector $\Pi_c\colon \widetilde{\cal SF}_c \to \big({\cal
  P}^0(c)\big)^d$:
\begin{align*}
  a(\bv - \Pi_c(\bv),\,\bv^0) = 0
  \qquad \forall \bv^0 \in \big({\cal P}^0(c)\big)^d,
\end{align*}
where $a(\bv,\bu) = (\K_c^{-1} \bv,\bu)_{L^2(c)}$.
The projector is computable using only the degrees of freedom \cite{Brezzi:2014}.
We define the consistency term as
$$
  a_{h,c}^{(1)}(\bu_h,\, \bv_h) = a(\Pi_c(\bu_h),\, \Pi_c(\bv_h))
  = \Int_c \K^{-1}_c \Pi_c(\bu_h) \cdot \Pi_c(\bv_h) \dx,
$$
and the stability term as
$$
  a_{h,c}^{(2)}(\bu_h,\, \bv_h) = (\bu_h - \Pi_c(\bu_h), \bv_h - \Pi_c(\bv_h))_{L^2(c)}.
$$
To generate a family of schemes, we can replace the $L^2$ inner
product with any other spectrally equivalent bilinear form.
Like in the mimetic framework, the consistency term is the exactness
property that holds when the first argument in the bilinear form is
the constant function:
$$
  a_{h,c}^{(1)}(\bu_h,\, \bv_h) 
  = \Int_c \K^{-1}_c \Pi_c(\bu_h) \cdot \Pi_c(\bv_h) \dx
  = \Int_c \K^{-1}_c \Pi_c(\bu_h) \cdot \bv_h \dx
  \qquad \forall \bv_h \in \widetilde{\cal SF}_h.
$$
Hence, the contibution of this term to the local mass matrix is like
in the mimetic scheme, $\R_c\,(\R_c^T\,\N_c)^{-1}\,\R^T_c$.

The conventional FE hybridization procedure works for the VEM as well.
The first equation in the variational formulation is replaced by a set
of cell-based equations:
\begin{align}
  a_{h,c}(\bu_{h,c},\, \bv_{h,c}) 
  - (\div_{h,c} \bv_{h,c},\,p_h)_{L^2(c)} 
  + \left<\bv_{h,c} \cdot \bn,\lambda_h\right>_{\partial c} = 0,
  \label{eq:VEM:local:1}
\end{align}
where $\bu_{h,c},\, \bv_{h,c} \in \widetilde{\cal SF}_c$,
$\div_c$ is the restriction of $\div_h$ to cell $c$, and 
$\lambda_h$ is the Lagrange multiplier.
%
It also holds that $\div_{h,c}\,\bv_{h,c} = \DIV_c\,\bv_c$ and
$$
  \left<\bv_{h,c} \cdot \bn,\lambda_h\right>_{\partial c} 
  = \Sum\limits_{f \in \partial c} \sigma_{c,f}\, |f|\, v_f^c\,\lambda_f,
$$
where $\lambda_f$ is the constant value of the Lagrange multiplier on
face $f$ (compare with~\bref{interface}).
The continuity equation is algebraically equivalent to \bref{continuity}.
Introducing vector $\blambda_c = \{\lambda_f\}_{f \in \partial c}$ and
recalling that $p_h$ has the constant value $p_c$ over cell $c$, the
hybridized equation can be rewritten as follows:
$$
  (p_c \bone - \blambda_c)^T\,\Sigma_c\,\F_c\,\bv_c
  = a_{h,c}^{(1)}(\bu_h,\, \bv_h),
$$
which also provides the relation $\bu_c = {\cal L}(p_c,\blambda_c)$.
The VEM is the linearity-preserving method.
According to Lemma~\ref{lem:equivalence}, the lowest-order
mixed-hybrid virtual element scheme is a member of the mimetic family
of schemes.

\section{Recent developments of the mimetic framework}
\label{sec:recent:developments}
\setcounter{equation}{0}

We highlight the new developments that extend the value of the mimetic framework 
for various physical applications.

\subsection{High-order schemes}

The development of a high-order mimetic scheme follows the same three
steps described above for the low-order mimetic schemes.
In the first step, we select the degrees of freedom that are
convenient for the definition of the primary divergence operator,
still denoted by $\DIV$.
With a slight abuse of notation we still use the symbols $\cC_h$ and
$\cF_h$ for the discrete spaces of pressure and flux unknowns,
respectively.

\begin{figure}[!t]
  \centering
  \includegraphics[scale=0.75]{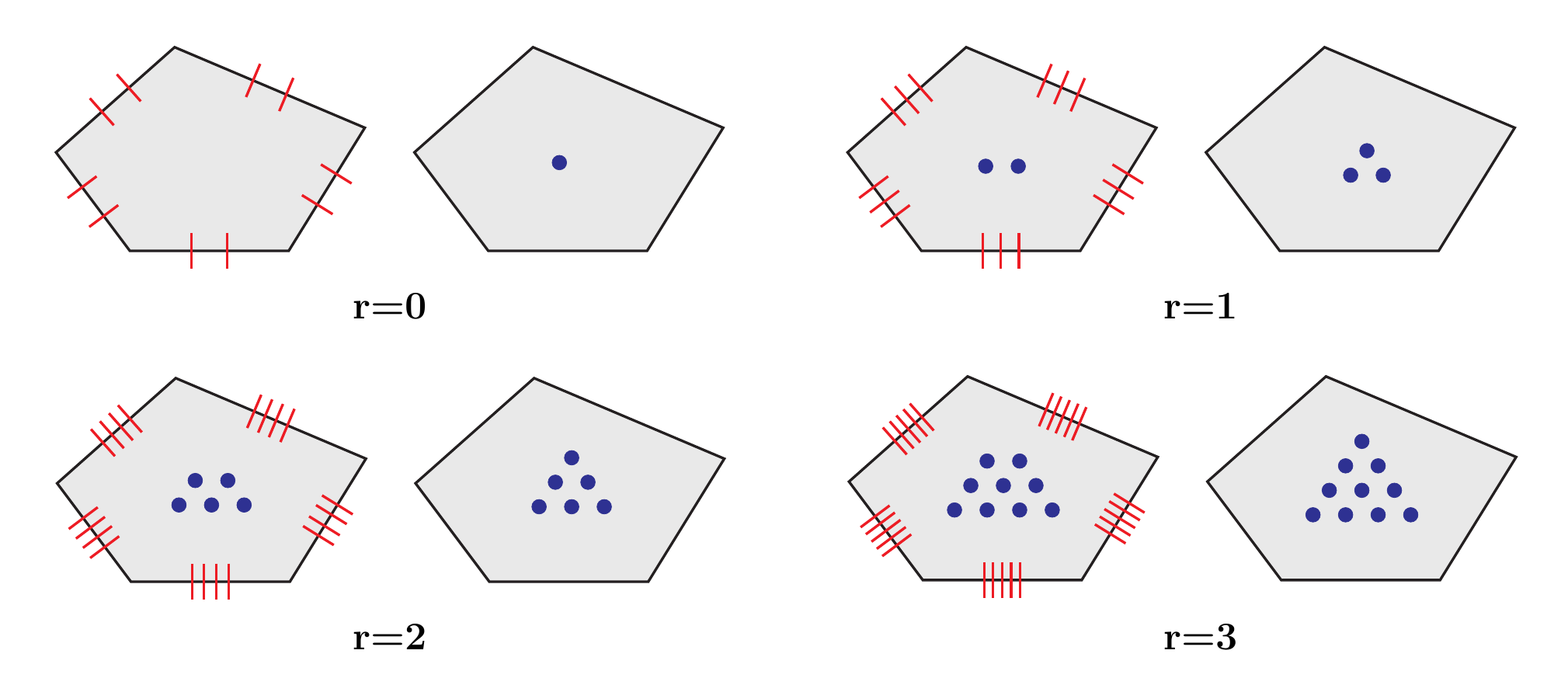}
  \caption{Degrees of freedom for $0\leq r\leq 3$ on a polygonal cell;
    for each polynomial degree $r$ we show the flux degrees of freedom
    on the left and the scalar degrees of freedom on the right.
    The edge/face moments of the normal component of the flux are
    denoted by a vertical line; the cell moments are denoted by a
    bullet.}
  \label{fig:dofs-HO}
\end{figure}

The discrete space $\cC_h$ contains multiple pressure unknowns that
can be associated with the solution moments up to order $r$.
The discrete space $\cF_h$ contains multiple flux unknowns both
associated with the mesh cells and the mesh faces.
The cell-based degrees of freedom represent moments of the flux up
to order $r$ except for the zero-th order moment.
The face-based degrees of freedom represent flux moments up to order
$r+1$, see Fig.~\ref{fig:dofs-HO}.

The discrete divergence operator $\DIV\colon \cF_h \to \cC_h$ is
defined cell-wise from the commutation property:
$$
  \big(\DIV\, \bu^I \big)_c = \DIV_c\, \bu_c^I = \big(\div\, \bu\big)^I_c,
$$
which is also a useful property for the error analysis.
This is the new design principle that can be generalized to other mimetic
operators.
The right-hand side is computable using only the degrees of freedom of
$\bu^I$.
Let $\psi \in {\cal P}^r(c)$ be a polynomial of order at most $r$.
Then, the definition of the moment and integration by parts give
$$
\big(\div\, \bu\big)^I_c
= \Frac{1}{|c|} \Int_c (\div\, \bu)\, \psi \dx
= -\Int_c \bu \cdot \nabla \psi\, \dx 
+ \Sum\limits_{f \in \partial c} \Int_f (\bu \cdot \bn_{c,f})\, \psi\, \dx.
$$

In the second step, we define the mimetic inner products in spaces
$\cC_h$ and $\cF_h$ as accurate approximations of the $L^2$ inner
products of pressure and flux functions.
The derivation is based on two high-order consistency conditions.
Since the local space $\cC_{h,c}$ is isomorphic to $\cP^r(c)$, the first
consistency condition is the obvious generalization
of~\eqref{eq:inner:product:scalar:low:order}:
$$
\big[p_c^I,\, q_c^I\big]_{\cC_{h,c}} = \int_c p\, q \dx
\qquad\forall p\in \cP^r(c),\ \forall q \in \cP^r(c).
$$
The consistency condition in space $\cF_{h,c}$ is defined as the following
exactness property:
$$
  \big[\bu_c^I,\, \bv_c^I\big]_{\cF_{h,c}} = \int_c \K_c^{-1} \bu\cdot \bv \dx
  \qquad\forall \bu\in \K_c\,\nabla\cP^{r+2}(c),\ \forall \bv \in {\cal SF}_c,
$$
where ${\cal SF}_c$ is a specially designed space containing the
vector functions $(\cP^{r+1}(c))^{d}$:
$$
{\cal SF}_c = \left\{ \bv\colon\ 
  \div\,\bv \in \cP^r(c),\quad 
  \bv\cdot \bn_f \in \cP^{r+1}(f)\quad \forall f \in \partial c \right\}.
$$

It is easy to show that the right-hand side of the consistency
condition is computable using the degrees of freedom introduced above.
Let $q\in \cP^{r+2}(c)$ be such that $\bu = \K_c\,\nabla q$. Then,
$$
  \int_c \K_c^{-1} \bu\cdot \bv \dx 
  = -\Int_c (\div \bv)\, q \dx 
    + \sum\limits_{f \in \partial c}\Int_f (\bv\cdot\bn_{c,f})\, q \dx.
$$
As $\bv$ is in ${\cal SF}_c$, the arguments of all the integrals in
the right-hand side above are polynomials.
Using the degrees of freedom of $\bv$ it is possible to reconstruct
$\div\,\bv$ inside $c$ and $\bv\cdot\bn_f$ on each $f\in\partial
c$, see~\cite{Gyrya-Lipnikov-Manzini:2015}, and all these integrals are
computable.
Combining the last formulas, we obtain the algebraic form of the
consistency condition:
$$
  \big[(\K_c \nabla q)^I,\, \bv_c^I\big]_{\cF_{h,c}} 
  = \big((\K_c \nabla q)^I\big)^T\, \M_{\cF,c}\, \bv_c^I
  = \big(\br_c(q)\big)^T\bv_c^I.
$$
To find a symmetric positive definite matrix $\M_{\cF,c}$, we need the analog of
Lemma~\ref{lem:NtR}, which obviously holds since the function $\K_c
\nabla \tilde q$ is in the space ${\cal SF}_c$ for any polynomial
$\tilde q \in \cP^{r+2}(c)$.

In the third step, we formulate the duality formula for the derived
gradient operator and use it in the global mimetic formulation.
The following error estimates have been shown in \cite{Gyrya-Lipnikov-Manzini:2015}:
$$
|||p^I - p_h |||_{\cC_h} + ||| \bu^I - \bu_h |||_{\cF_h} \le C\, h^{r+2}.
$$

\begin{remark}
  In the case when the diffusion tensor $\K$ is no longer constant,
  the consistency condition has to be modified.  Let $\Pi_c^r$ denote
  the local $L^2$ projector on the space of polynomial functions of
  order $r$.  Then, the modified consistency condition reads:
  $$
  \big[\big(\Pi_c^{r+1}(\K \nabla q)\big)^I,\, \bv_c^I\big]_{\cF_{h,c}} = \int_c \nabla q \cdot \bv \dx
  \qquad\forall q\in \cP^{r+2}(c),\ \forall \bv \in {\cal SF}_c.
  $$
  After that, the inner product matrix $\M_{\cF,c}$ is derived
  following the same steps.
\end{remark}

\subsection{Nonlinear parabolic problems}
\label{sec:divk}

The consistency term in the formula for matrix $\M_{\cF,c}$ (see
\bref{consistency}) contains the inverse of the diffusion tensor.
Therefore, numerical difficulties may arise in solving nonlinear
parabolic problems of type
$$
  \frac{\partial p}{\partial t} - \div(k(p)\, \nabla p) = b,
$$
where function $k(p)$ cannot be uniformly bounded from below.  
For instance, on a uniform one dimensional mesh, the numerical flux at
mesh point $x_i$ is proportional to the difference of the neighboring
pressures and the transmissibility coefficient $T_i$:
$$
  u_i = -T_i \Frac{p_{i+1/2} - p_{i-1/2}}{h},
  \qquad
  T_i = \Frac{2\,k_{i-1/2}\, k_{i+1/2}}{k_{i-1/2} + k_{i+1/2}}.
$$
If $k_{i-1/2} \ll k_{i+1/2}$, the numerical flux goes to zero as
$k_{i-1/2} \to 0$ and may lead to a nonphysical solution as shown by
the numerical experiment considered in
Section~\ref{eq:Marshak:heat:equation} (see also \cite{Maslyankin:1977}).
To obtain an accurate solution, we have to replace the harmonic
average with the arithmetic average in the definition of the
transmissibility coefficient.
A possible strategy in the mixed finite element framework (see, e.g.
\cite{Arbogast:1997}) consists in using two velocity variables, $\bv
= -\nabla p$ and $\bu = k(p) \bv$.  
However, the corresponding weak formulation cannot have 
face-based equations of type $v_c^f = k_f u_c^f$ which are natural from a physical viewpoint.
Here $k_f$ is the face-based diffusion coefficient.
Such face-based equations lead to an algebraic problem that can be symmetrized
only on special meshes.
We describe a new mimetic scheme that allows us to use different
values $k_f$ on different mesh faces.
The mimetic framework always guarantees the symmetry of the resulting
algebraic problem.


\subsubsection{A new pair of primary and derived mimetic operators}
Let us consider a more general form of the diffusion coefficient, $\K\, \k(p)$, 
where $\K$ is a discontinuous tensor independent of $p$ and $\k(p)$ is
a discontinuous scalar function of $p$.
The underlying mixed formulation is
\begin{equation}\label{parabolic}
  \begin{array}{rcl}
    \bu &=& -(\K\,\nabla) p,\\[0.1ex]
    \Frac{\partial p}{\partial t} + (\div\,\k) \, \bu &=& b.
  \end{array}
\end{equation}

The combined operators $\div\,\k$ and $\K\,\nabla$ are dual to each
other with respect to the weighted $L^2$ inner products by using
$\k\K^{-1}$ as weight:
\begin{equation}
  \label{Green:duality}
  \Int_\Omega (\div\,\k \bu)\, q \dx
  = -\Int_\Omega \k\,\K^{-1} \bu \cdot (\K\,\nabla)\, q \dx
  \qquad \forall \bu \in H_{div}(\Omega),\  \forall q \in H^1_0(\Omega).
\end{equation}

Consider again the three-step construction of the mimetic framework.
In the \emph{first} step we need to specify the degrees of freedom.
For the pressure variable, we consider the same discrete space $\cC_h$
of grid functions that consist of one value per cell and the same
discrete space $\Lambda_h$ of grid functions that consist of one value
per face.
The discrete space $\cF_{h}$ has the same dimension as in the linear
case, but the discrete fluxes in $\cF_{h}$ obey a different continuity
condition:
\begin{equation}\label{k-continuity}
  \k_f^{c_1}\, u^{c_1}_f = \k_f^{c_2}\, u^{c_2}_f
\end{equation}
on each interior face $f$ shared by cells $c_1$ and $c_2$.
Here, $\k_f^{c_1}$ and $k_f^{c_2}$ are accurate one-side
approximations of the diffusion coefficient $\k$.
For example, in regions where function $k$ is continuous, we can take 
$\k_f^{c_1} = \k_f^{c_2} = \k_f$ as a weighted average of the cell-centered 
values $\k(p_{c_1})$ and $\k(p_{c_2})$ calculated using the most recent 
approximation to solution $p$ in the cells $c_1$ and $c_2$, respectively.
The weights are the distances between the cell centers and face $f$.

The primary mimetic operator approximates the combined operator
$\div\,\k$.
It is defined locally on each mesh cell 
using a straightforward discretization of the divergence theorem 
(compare with formula \bref{primary-div}):
\begin{equation}\label{primary-divk}
  \left.(\DIV^\k\, \bu_h)\right|_c 
  \equiv \DIV^\k_c\, \bu_c 
  = \frac{1}{|c|} \Sum\limits_{f \in \partial c} \sigma_{c,f}\,|f|\, \k_f^c u^c_f.
\end{equation}

Since $\bu_h$ is an algebraic vector, it is convenient to think about the discrete
divergence operator $\DIV^\k\colon \cF_h \to \cP_h$ as a matrix acting 
between two spaces.
This matrix has full rank when $\k_f^c > 0$.

Let us introduce a cell-based diagonal matrix ${\cal K}_c$ formed by coefficients $\k_f^c$,
$f \in \partial c$.
Then, the primary mimetic operators in \bref{primary-div} and \bref{primary-divk}
can be connected as follows:
$$
  \DIV^\k_c\, \bu_c 
  = (\DIV_c\, {\cal K}_c) \bu_c.
$$

The {\it second} step in the mimetic discretization framework is to define
the inner products in spaces $\cC_h$ and $\cF_h$ that are accurate approximations 
of the integrals in \bref{Green:duality}.
Such inner products can be defined again cell-by-cell.
Moreover, the inner product in space $\cC_h$ can be defined as in
Section~\ref{subsec:mimetic:discretization:framework} by the relation
$\big[p_c,\, q_c\big]_{\cC_{h,c}} = |c|\, p_c\, q_c$.

The weight in the other $L^2$ inner product is given by $\k\,\K^{-1}$.
By our assumption, $\K$ is a piecewise constant tensor on mesh $\Omega_h$.
Instead, the scalar coefficient $\k$ can be a quite general
non-negative function.
An acceptable first-order error is committed when we replace $\k$ by
the piecewise constant function with value $\k_c = \k(p_c)$ in cell
$c$.

The consistency condition is expressed through the exactness property:
$$
  \big[\bu_c^I,\, \bv_c^I\big]_{\cF_{h,c}} 
   = \int_c \k(p_c)\, \K_c^{-1} \bu\cdot \bv \dx
  \qquad\forall \bu\in (\cP^0(c))^{d},\ \forall \bv \in {\cal SF}_c.
$$
Since $\k_c\, \K^{-1}_c$ is a constant tensor in cell $c$, the
derivation of the inner product matrix (still denoted by $\M_{\cF,c}$)
proceeds as in Section~\ref{subsec:mimetic:discretization:framework}.

The {\it third} step in the construction of the mimetic framework is
to obtain the formula for the derived operator, which, in this case,
is an approximation of the combined operator $\K\,\nabla$.
The continuum Green formula for cell $c$ is given by
\begin{equation}
  \label{IBP_local_general}
  \Int_c (\div\k\bu)q\, dx
  - \Int_{\p c} (\k\bu\cdot\bn)\,q\, dx
  =
  - \Int_c\k\,\K^{-1}(\K\,\nabla q)\cdot\bu\,dx.
\end{equation}
Now, the derived operator $\GRADt_c\colon \cC_{h,c} \times
\Lambda_{h,c} \to \cF_{h,c}$ satisfies the discrete integration by
parts formula
\begin{align*}
  \big[\DIV^\k_c\, \bu_c,\, q_c\big]_{\cC_{h,c}} 
  - \Sum\limits_{f\in\p c}\sigma_{c,f}\,|f|\,k_f^c\, u^c_f \lambda_f
  = -\Big[\bu_c,\, \GRADt_c\,
  \Big(\!\!\begin{array}{c}q_c \\ \blambda_c \end{array}\!\!\Big)
  \Big]_{\cF_{h,c}}
\end{align*}
for all $\bu_c\in \cF_{h,c}$,  $q_c \in \cC_{h,c}$, and $\blambda_c \in \Lambda_{h,c}$.
This local mimetic formulation gives the following formula for the physical fluxes:
\begin{equation}\label{fluxes}
  {\cal K}_c 
  \left(\!\!\!\begin{array}{c}u^c_{f_1} \\ \vdots \\ u^c_{f_{n_c}}\end{array}\!\!\!\right)
  =
  -{\cal K}_c\,\GRADt_c\, \Big(\!\!\begin{array}{c}p_c \\ \blambda_c \end{array}\!\!\Big) 
  =
  {\cal K}_c \M_{\cF,c}^{-1} {\cal K}_c 
    \left(\!\begin{array}{c} 
    \sigma_{c,f_1}\, |f_1|(p_c - \lambda_{f_1}) \\ 
    \vdots \\
    \sigma_{c,f_{n_c}}\, |f_{n_c}|(p_c - \lambda_{f_{n_c}})  \end{array}\!\right).
\end{equation}
Note that this formula uses the symmetric matrix ${\cal K}_c \M_{\cF,c}^{-1} {\cal K}_c$.
Our numerical experiments show that the resulting scheme is second-order accurate.

\subsubsection{Marshak heat equation}
\label{eq:Marshak:heat:equation}
Let us consider the modified Marshak heat equation
\cite{Samarskii:1963,Maslyankin:1977} in the rectangular domain
$(0,3)\times(0,1)$ with zero source term, $\K = \I$, and $\k(p) =
p^3$.
The initial value is $p(\bx,0) = 10^{-3}$.
We set the time-dependent Dirichlet boundary condition
$p=p^0(0,t)=0.78\,t^{1\slash{3}}$ on the left side of $\Omega$, the
constant boundary condition $p=10^{-3}$ on the right side, and the
homogeneous Neumann boundary conditions on the remaining sides.

\begin{figure}[h!]
  \centering
  \includegraphics[scale=0.75]{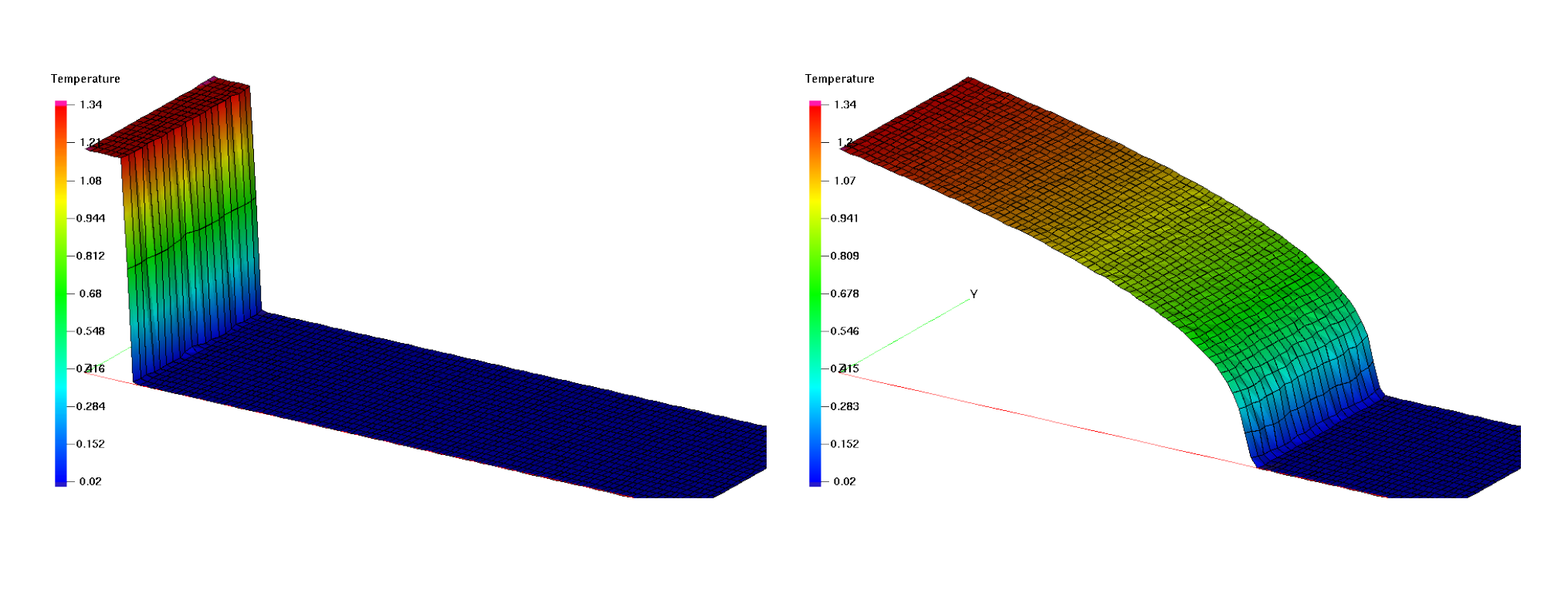}
  \caption{Solution snapshots at time $t=5.0$ for the standard (left)
    and new (right) MFD schemes. The solution in the right panel shows
    the correct and accurate position of the wave front at the chosen
    time.}
  \label{fig:marshak-new}
\end{figure}

We solve the parabolic equation on the randomly perturbed
quadrilateral mesh that have three times more cells in the x-direction
than in the y-direction.
We use the backward Euler time integration scheme and the weighted
arithmetic average definition of the face-based diffusion coefficients
$\k_f$ described above.
To reduce the impact of the time integration error, we use small time
steps.
Comparison of pictures in Fig.~\ref{fig:marshak-new} show that the new
scheme fixes the deficiencies of the old MFD scheme and leads to the
correct speed of propagation of the non-linear wave.

\section{Conclusions}
\label{sec:conclusions}

We described a few design principles used in the derivation of mimetic schemes 
for the numerical solution of PDEs.
We established the bridges with a few FV and FE methods by showing
that three popular discretization frameworks (MFD, FV and FE) use the equivalent 
design principles which leads to algebraically equivalent schemes.
We illustrated the flexibility of the mimetic discretization framework 
to tackle challenging numerical issues in computer modeling of engineering
problems with two examples: derivation of higher-order schemes and convergent 
schemes for nonlinear problems with small diffusion coefficients.

\section*{Acknowledgments}

This work was carried out under the auspices of the National Nuclear
Security Administration of the U.S. Department of Energy at Los Alamos
National Laboratory under Contract No. DE-AC52-06NA25396.  
The authors acknowledge the support of the US Department of Energy
Office of Science Advanced Scientific Computing Research (ASCR)
Program in Applied Mathematics Research.

The meshes for the Marshak problem were created and managed using the
mesh generation toolset MSTK (software.lanl.gov/MeshTools/trac)
developed by Dr. Rao Garimella at Los Alamos National Laboratory.


\end{document}